\documentclass[12pt,a4paper]{amsart}

\usepackage{graphicx}
\usepackage{amsfonts,amsmath,latexsym,amssymb,amsthm}

\newcounter{theoremcounter}
\newcounter{lemmacounter}

\newcounter{dummycounter}

\newcounter{corcounter}

\newcounter{emptycounter}
\newcounter{defcounter}

\newtheorem{theorem}[theoremcounter]{Theorem}

\newtheorem{lemma}[lemmacounter]{Lemma}

\newtheorem{korollar}[corcounter]{Corollary}
\newtheorem{remark}{Remark}
\newtheorem{definition}[defcounter]{Definition}

\numberwithin{equation}{section}
\numberwithin{lemmacounter}{section}
\numberwithin{propcounter}{section}
\numberwithin{corcounter}{section}
\numberwithin{conjcounter}{section}
\numberwithin{theoremcounter}{section}
\numberwithin{probcounter}{section}

\newcounter{eqncounter}

\numberwithin{equation}{eqncounter}

\newcommand\pw{\mathfrak{p}}
\newcommand\f{\iota}

\newcommand\coll{\mathcal{C}}
\newcommand\collg{\coll_e^{(g)}}

\newcommand\Mainterm{D_{\en_K}}
\newcommand\Errorterm{B_{\en_K}}
\newcommand\IG{F}

\newcommand\ALS{ALS}
\newcommand\m{m}
\newcommand\wK{w_K}
\newcommand\M{M}

\newcommand\IR{\mathbb R}
\newcommand\IC{\mathbb C}
\newcommand\IZ{\mathbb Z}

\newcommand\IP{\mathbb P}
\newcommand\IQ{\mathbb Q}

\newcommand\en{\mathcal{N}}
\newcommand\hen{H_\mathcal{N}}
\newcommand\balf{{\mbox{\boldmath $\alpha$}}}
\newcommand\vnull{{\mbox{\boldmath $0$}}}

\newenvironment{rproof}{\addvspace{\medskipamount}\par\noindent{\it Proof.\/}}
{\unskip\nobreak\hfill$\Box$\par\addvspace{\medskipamount}}
\newcommand\ord{\mathop{\rm ord}\nolimits}

\newcommand\Oseen{{\mathcal{O}}}
\newcommand\A{{\mathfrak{A}}}
\newcommand\C{{\mathfrak{C}}}
\newcommand\D{{\mathfrak{D}}}

\renewcommand{\L}{{\mathfrak{L}}}

\newcommand\Lamen{\Lambda_{\en}}
\newcommand\LamenK{\Lambda_{\en_K}}
\newcommand\Qbar{\overline{\IQ}}

\newcommand\kbar{\overline{k}}

\newcommand\vx{{\bf x}}
\newcommand\vy{{\bf y}}
\newcommand\vz{{\bf z}}

\newcommand\henK{H_{{\mathcal{N}_{K}}}}
\newcommand\Cl{Cl}
\newcommand\Ce{D}
\newcommand\Da{D}

\newcommand\cS{c}
\newcommand\CS{C}
\newcommand\cDe{c_{\Delta}}
\newcommand\ceins{c_1}
\newcommand\czwei{c_2}

\begin{document}\baselineskip=17pt
\title{Counting points of fixed degree and bounded height}

\author{Martin Widmer}

\address{Department of Mathematics\\ 
The University of Texas at Austin\\ 
1 University Station C1200\\
Austin, Texas 78712\\ 
U.S.A.}

\email{widmer@math.utexas.edu}

\date{May 15, 2009}

\subjclass[2000]{Primary 11R04; Secondary 11G50, 11G35}

\keywords{Height, Northcott's Theorem, counting}

\begin{abstract}
We consider the set of points in projective $n$-space that generate an extension of degree $e$
over given number field $k$, and deduce an asymptotic formula for the number of such points of absolute
height at most $X$, as $X$ tends to infinity. We deduce a similar such formula with instead of the 
absolute height, a so-called adelic-Lipschitz height.
\end{abstract}

\maketitle

\section{Introduction}
Let $k$ be a number field of degree $\m=[k:\IQ]$ in a fixed algebraic closure $\kbar$ of $k$
and let $n$ be a positive integer. Write $\IP^n(\kbar)$ for the projective space of dimension $n$ over the field $\kbar$ and
denote by $H(\cdot)$ the non-logarithmic absolute Weil height on $\IP^n(\kbar)$ as defined in \cite{BG} p.16.
A fundamental property of the height, usually associated with the name of Northcott due to his result Theorem 1 in \cite{26},
states that subsets of $\IP^n(\kbar)$ of bounded degree and bounded height are finite.
This raises the question of estimating the cardinality of such a set as the height bound gets large. 
Schanuel proved in \cite{25} that for the counting function
\begin{alignat*}1
Z_H(\IP^n(k),X)=|\{P\in \IP^n(k);H(P)\leq X\}|
\end{alignat*}
one has an asymptotic formula
\begin{alignat}1
\label{Scha1}
Z_H(\IP^n(k),X)=S_k(n)X^{\m(n+1)}+O(X^{\m(n+1)-1}\log X)
\end{alignat}
as $X$ tends to infinity
where $S_k(n)$ is a positive constant involving all classical field 
invariants (see (\ref{Schanuelconst}) for its definition) and the constant implied 
by the Landau $O$-symbol depends on $k$ and $n$.
The logarithm can be omitted in all cases except for
$n=\m=1$.\\

A projective point $P=(\alpha_0:...:\alpha_n)$ in $\IP^n(\kbar)$ has a natural degree
defined as 
\begin{alignat*}1
[k(P):k]
\end{alignat*}
where $k(P)$ denotes the
extension we get by adjoining all ratios $\alpha_i/\alpha_j$
$(0\leq i,j\leq n, \alpha_j\neq 0)$ to $k$.
In 1993 Schmidt drew attention to the more general set 
\begin{alignat*}1
\IP^n(k;e)=\{P\in \IP^n(\kbar);[k(P):k]=e\}
\end{alignat*} 
of points with relative degree $e$.
Clearly $\IP^n(k;1)=\IP^n(k)$ and so Schanuel deals with the case $e=1$.
For the counting function
\begin{alignat*}1
Z_H(\IP^n(k;e),X)=|\{P\in \IP^n(k;e);H(P)\leq X\}|
\end{alignat*}
Schmidt \cite{22} proved the following general estimates
\begin{alignat}1
\label{ThSchm1'}
\cS X^{\m e(\max\{e,n\}+1)}\leq Z_H(\IP^n(k;e),X)&\leq \CS X^{\m e(e+n)}
\end{alignat}
where $c=c(k,e,n)$ and $C=C(k,e,n)$ are positive constants depending solely
on $k,e$ and $n$.
The upper bound holds for $X\geq 0$ and the lower
bound holds for $X\geq X_0(k,e,n)$ depending also on $k,e$ and $n$.
Moreover one can choose $C=2^{\m e(e+n+3)+e^2+n^2+10e+10n}$.\\

For $k=\IQ$ more is known. 
Schmidt \cite{14} investigated the quadratic case. 
Here he provided not only the correct order of magnitude but he found also the 
precise asymptotics and this in all dimensions $n$.
As $X$ tends to infinity one has
\begin{alignat}3\label{ThSchm2}
Z_H(\IP^n(\IQ;2),X)=
\left\{ \begin{array}{llc}
\Ce_1X^{6}+O(X^4\log X)& \mbox{if $n=1$} \\
\Ce_2X^{6}\log X+O(X^6\sqrt{\log X})& 
\mbox{if $n=2$}\\
\Ce_nX^{2(n+1)}+O(X^{2n+1})& \mbox{if $n>2$}
\end{array}\right..
\end{alignat}
The constant implied by the $O$-symbol depends only on $n$.
In fact Schmidt's result was more precise since it gave the
asymptotics for real and imaginary quadratic points separately.
Here $\Ce_1=\frac{8}{\zeta(3)}$, 
$\Ce_2=\frac{8(12+\pi^2)}{\zeta(3)^2}$ and $\Ce_n=\Ce(\IQ,2,n)$ is given by
\begin{alignat*}1
\Ce(\IQ,2,n)=\sum_K S_K(n) 
\end{alignat*}
where the sum runs over 
all quadratic fields $K$.
Schmidt proved also an analogue to the above result for a more general kind
of height and showed that this leads to asymptotic formulae for
the number of decomposable quadratic forms (i.e. product of two linear forms) 
$f(x_0,...,x_n)=\sum_{0\leq i\leq j\leq n}a_{ij}x_ix_j$
with coefficients $a_{ij}$ in $\IZ$ having $|a_{ij}|\leq X$ and moreover for the number of symmetric $(n+1)\times(n+1)$
matrices with rank $\leq 2$ such that $b_{ii}\in \IZ$, $|b_{ii}|\leq X$
and $2b_{ij}\in \IZ$, $2|b_{ij}|\leq X$ for $i\neq j$.
Already way back in 1967 Schmidt \cite{86} introduced more general classes
of heights where the maximum norm in (\ref{height}) at the infinite places is replaced by an arbitrary but fixed distance function.
More recently Thunder \cite{ThMiHl} and Roy-Thunder \cite{79} introduced ``twisted heights'' 
which allow also modifications at the finite places.\\

One year after Schmidt's article on quadratic points Gao \cite{7} made further progress. 
He proved that if $n>e>2$ one has 
\begin{alignat}1
\label{ThGao}
Z_H(\IP^n(\IQ;e),X)=\Ce(\IQ,e,n)X^{e(n+1)}+O(X^{e(n+1)-1})
\end{alignat}
as $X$ tends to infinity. The constant implied by the $O$-symbol depends on
$e$ and $n$ and the constant $\Ce(\IQ,e,n)$ is given by $\sum_K S_K(n)$ where the sum runs over 
all extensions $K$ of $\IQ$ of degree $e$.
For $1\leq n\leq e$ Gao showed that the correct order 
of magnitude of $Z_H(\IP^n(\IQ;e),X)$ is $X^{e(e+1)}$.
Here the asymptotics are still unknown, even in the case
$e=3$ and $n=2$ of cubic points in two dimensions.\\
 
Schmidt's and Gao's results are restricted to the ground field $k=\IQ$. 
A more recent result of Masser and Vaaler \cite{1}
gives the asymptotics for the number of points of fixed degree over an arbitrary
fixed number field $k$,
but only in dimension $n=1$. Masser and Vaaler established the asymptotic formula
\begin{alignat}1\label{ThMV1}
Z_H(\IP^1(k;e),X)=
eV_{\IR}(e)^{r_k}V_{\IC}(e)^{s_k}S_k(e)X^{\m e(e+1)}+O(X^{\m e(e+1)-e}\log X)
\end{alignat}
as $X$ tends to infinity.
The constants $V_{\IR}(e), V_{\IC}(e)$ have their origins
in \cite{43}. Moreover the logarithm   
can be omitted in all cases except $(\m,e)=(1,1)$ and
$(\m,e)=(1,2)$ and the constant implied by the $O$-symbol depends on
$k$ and $e$. Unfortunately the proof of Masser and Vaaler's
theorem shed no light on the case $n>1$.
Very roughly speaking Masser and Vaaler's idea was to interpret the height of the root of an irreducible polynomial in $k[x]$ of fixed degree
$e$ as a suitable height of the coefficient vector
of this polynomial and to proceed by counting minimal
polynomials with respect to this modified height. 
To carry out this plan they had to generalize the
class of heights introduced by Schmidt \cite{86}, allowing now
different distance functions at the infinite places
instead of only one for all infinite places as Schmidt did. 
On the other hand Masser and Vaaler
had to impose a technical condition, associated
with the name of Lipschitz, on the
boundaries of the unit balls given by the respective 
distance function. They therefore introduced so-called
Lipschitz systems, giving what one could call
Lipschitz heights.\\

In the present article we establish the
asymptotics for $\IP^n(k;e)$ if $n$ is slightly larger than $5e/2$.
Let us write
\begin{alignat*}1
\Ce=\Ce(k,e,n)=\sum_K S_K(n)
\end{alignat*} 
for the formal sum taken over all extensions of $K$ of $k$ in $\kbar$ of degree $e$.
We have the following result.
\begin{theorem}\label{mainthintro}
Let $e,n$ be positive integers and $k$ a number field
of degree $\m$ and suppose that 
$n>5e/2+4+2/(\m e)$.
Then the sum defining $\Ce$ converges and 
as $X$ tends to infinity we have
\begin{alignat*}1
Z_{H}(\IP^n(k;e),X)=\Ce X^{\m e(n+1)}
+O(X^{\m e(n+1)-1}\log X).
\end{alignat*}
The logarithm can be omitted unless $(\m e,n)=(1,1)$
and the constant implied by the $O$-symbol depends on
$k,e$ and $n$.
\end{theorem}
If $e$ and $n$ are both larger than one there is a considerable gap between the exponents of the lower and
the upper bound in (\ref{ThSchm1'}). Schmidt mentioned that the lower bound is likely to be nearer the truth
than the upper bound. Our Theorem \ref{mainthintro} confirms Schmidt's conjecture at least if $n$ is large enough.
We will prove a more general result involving adelic-Lipschitz heights.\\

Let us give a single new example of our theorem.
We take $n=11$, $k=\IQ(i)$, $e=2$, so that we 
are counting points in eleven dimensions quadratic
over $\IQ(i)$. For the number 
$Z=Z_H(\IP^{11}(\IQ(i);2),X)$ of points of height at most $X$, the Schmidt bounds are 
$X^{48}\ll Z\ll X^{52}$
for $X\geq X_0$,
with absolute implied constants. Our result implies
that 
\begin{alignat*}1
Z=DX^{48}+O(X^{47})
\end{alignat*}
with
\begin{alignat*}1
D=12\cdot(2\pi)^{24}\sum_{K\atop [K:\IQ(i)]=2}\frac{h_KR_K}{\wK \zeta_K(12)|\Delta_K|^6}.
\end{alignat*}
\indent Our proof follows the general strategy of Schmidt and Gao.
Their audacious idea was
to prove a result similar to (\ref{Scha1})
but with $\IP^n(K)$ replaced by $\IP^n(K/\IQ)$ the subset
of primitive points in
$\IP^n(K)$; by definition these satisfy $K=\IQ(P)$. Now $\IP^n(\IQ;e)$ is a disjoint union
of the sets $\IP^n(K/\IQ)$ where $K$ runs over all
number fields of degree $e$. 
For each $\IP^n(K/\IQ)$ the
main term is the same as that in (\ref{Scha1}) with $K$ instead of $k$, but for $e=2$
Schmidt obtained a more precise error term
\begin{alignat}1\label{Schmerror1}
O\left(\frac{\sqrt{h_KR_K\log(3+h_KR_K)}}{|\Delta_K|^{n/2}}X^{2n+1}\right)
\end{alignat}
where the constant in $O$ depends only on $n$ but is independent
of the field $K$. This is the major step of the proof
and involves many new ideas.
Now one can sum over all quadratic number fields and the 
Theorem of Siegel-Brauer ensures that the
sum over the main terms $S_K(n)$ as well as over
the error terms converges provided $n>2$. For similar reasons the restriction $n>e$ in Gao's result appears.\\

We close the introduction with some remarks about the structure of the paper.
First we take up the definition
of an adelic-Lipschitz system from \cite{art1} on a number field
and we define a uniform adelic-Lipschitz system 
on the collection of all extensions of $k$ of degree $e$. 
This then gives
rise to a class of heights $\hen$ defined
on $\IP^n(k;e)$.
The main result asymptotically
estimates the counting function of $\IP^n(k;e)$
with respect to the height $\hen$.
In Theorem \ref{mainthintro} we used only the simplest formulation by
choosing a special uniform adelic-Lipschitz system
with maximum norms at all places so that the corresponding adelic-Lipschitz height $\hen$ becomes just the Weil height $H$.
In Section \ref{introchap4} we state our main theorem for general adelic-Lipschitz heights. It is in \cite{art3} and \cite{art4}
where we see the advantage of working in such generality.
In \cite{art3} we are concerned with counting points of fixed degree on linear subvarieties
of projective space. In \cite{art4} we prove the following result: 
let $m,n$ be positive integers with $n>\max\{6m+2+2/m,m^2+m\}$. Then
as $X$ tends to infinity the number 
of algebraic numbers $\alpha$ of degree $mn$ such that $\IQ(\alpha)$
contains a subfield of degree $m$ and $H(1,\alpha)\leq X$ is
asymptotically equal to
\begin{alignat*}1
\Ce'(m,n) X^{mn(n+1)}
\end{alignat*}
where $\Ce'(m,n)=\sum_K nV_{\IR}(n)^{r_k}V_{\IC}(n)^{s_k}S_k(n)$
and the sum runs over all number fields of degree $m$.
Note that the subfield condition reduces the order of magnitude
from $X^{mn(mn+1)}$ to $X^{mn(n+1)}$.\\
 
In Section \ref{2subsec5} we prove the main result Theorem \ref{main theorem} which is a
generalization of Theorem \ref{mainthintro} to adelic-Lipschitz heights.
Section \ref{countingfields} is devoted to some simple lower and upper bounds for the number of extensions
$K/k$ of fixed degree with $\delta(K/k)\leq T$, where $\delta(K/k)=\underset{\alpha}\inf\{H(1,\alpha);K=k(\alpha)\}$. The invariant $\delta(K/k)$ 
was already introduced in \cite{art1}. Our bounds are essentially by-products of the
proof of Theorem \ref{main theorem}.

\section*{Acknowledgements}
I would like to thank my Ph.D. advisor David Masser for his constant support.
His comments led to many significant improvements of this article.
I also would like to thank Jeffrey Vaaler for helpful discussions
and Gao Xia for showing us his Ph.D. thesis. 
This work was financially supported by the Swiss National Science Foundation.

\section{Adelic-Lipschitz heights revisited}\label{2subsec4}
The Subsections \ref{prelim1}, \ref{subsecdefALS} and \ref{2subsec2} of this section are entirely contained in \cite{art1}.
But adelic-Lipschitz heights are crucial for the entire paper and thus, for convenience of the reader,
we introduce this notion here once again.
Before we can define adelic-Lipschitz heights we have to fix some basic notation.
For a detailed account on heights we refer to \cite{BG} and \cite{3}.\\

\subsection{Preliminaries}\label{prelim1}
Let $K$ be a finite extension
of $\IQ$ of degree $d=[K:\IQ]$. By a place $v$ of $K$ we mean an equivalence class of non-trivial absolute values on $K$.
The set of all places of $K$ will be denoted by $M_K$.
For each $v$ in $M_K$ we write $K_v$ for the completion
of $K$ at the place $v$ and $d_v$
for the local degree defined by
$d_v=[K_v:\IQ_v]$ 
where $\IQ_v$ is a completion with respect to the place which extends to $v$. A place $v$ in $M_K$ corresponds either to a non-zero prime ideal $\pw_v$
in the ring of integers $\Oseen_K$ or to an embedding
$\sigma$ of $K$ into $\IC$.
If $v$ comes from a prime ideal we call $v$
a finite or non-archimedean place and denote this by $v\nmid \infty$ and if $v$ corresponds to
an embedding we say $v$ is an infinite or archimedean
place and denote this by $v\mid \infty$.
For each place in $M_K$ we choose a representative
$|\cdot|_v$,
normalized in the following way:
if $v$ is finite and $\alpha\neq 0$ we set by convention
\begin{alignat*}3
|\alpha|_{v}=N\pw_v^{-\frac{\ord_{\pw_v}(\alpha\Oseen_K)}{d_v}}
\end{alignat*}
where $N\pw_v$ denotes the norm of $\pw_v$
from $K$ to $\IQ$ and $\ord_{\pw_v}(\alpha\Oseen_K)$
is the power of $\pw_v$ in the prime ideal decomposition
of the fractional ideal $\alpha\Oseen_K$.
Moreover we set
\begin{alignat*}3
|0|_{v}=0.
\end{alignat*}
And if $v$ is infinite and corresponds to an embedding $\sigma:K \hookrightarrow \IC$ we define
\begin{alignat*}3
|\alpha|_{v}=|\sigma(\alpha)|.
\end{alignat*}
If $\alpha$ is in $K^*=K\backslash\{0\}$ then 
$|\alpha|_v\neq 1$
holds for only a finite number of places $v$.\\

Throughout this article $n$ will denote a positive rational integer. The height
on $K^{n+1}$ is defined by
\begin{alignat}3
\label{height}
H(\alpha_0,...,\alpha_n)=\prod_{M_K}\max\{|\alpha_0|_v,...,|\alpha_n|_v\}^{\frac{d_v}{d}}.
\end{alignat}
Due to the remark above this is in fact a finite product. 
Furthermore this definition is independent of the field $K$ containing
the coordinates (see \cite{BG} Lemma 1.5.2 or \cite{3} pp.51-52) and therefore
defines a height on $\Qbar^{n+1}$ for an algebraic closure
$\Qbar$ of $\IQ$.
The well-known \it product formula \em (see \cite{BG} Proposition 1.4.4) asserts that
\begin{alignat*}3
\prod_{M_K}|\alpha|_v^{d_v}=1 \text{ for each $\alpha$ in $K^*$}.
\end{alignat*}
This implies in particular that the value of the height in (\ref{height}) does not change if we multiply each coordinate with a fixed element of $K^*$. 
Therefore one can define a height on points
$P=(\alpha_0:...:\alpha_n)$ in $\IP^n(\Qbar)$ by
\begin{alignat}3
\label{heightproj}
H(P)=H(\alpha_0,...,\alpha_n).
\end{alignat}
Moreover, to evaluate the height, we can assume that one of the coordinates is $1$ which shows that $H(\balf)\geq 1$ for $\balf\in \Qbar^{n+1}\backslash\{\vnull\}$.
The equations (\ref{height}) and (\ref{heightproj}) define the absolute non-logarithmic
projective Weil height or just Weil height.\\

\subsection{Adelic-Lipschitz systems on a number field}\label{subsecdefALS}
Let $r$ be the number of  
real embeddings and $s$ the number of pairs of complex conjugate embeddings
of $K$ so that $d=r+2s$.
Recall that $M_K$ denotes the set of places of $K$.
For every place $v$ we fix a
completion $K_v$ of $K$ at $v$. The value set of $v$, $\Gamma_v:=\{|\alpha|_v;\alpha \in K_v\}$
is equal to $[0,\infty)$ if $v$ is archimedean,
and to
\begin{alignat*}3
\{0,(N\pw_v)^{0},(N\pw_v)^{\pm 1/d_v},(N\pw_v)^{\pm 2/d_v},...\}
\end{alignat*}
if $v$ is non-archimedean.
For $v \mid \infty$ we identify $K_v$ with $\IR$ or
$\IC$ respectively and we identify $\IC$ with
$\IR^2$ via $\xi\longrightarrow (\Re(\xi),\Im(\xi))$
where we used $\Re$ for the real and $\Im$ for the
imaginary part of a complex number.\\

Before we can introduce adelic-Lipschitz systems we have to 
give a technical definition.
For a vector $\vx$ in $\IR^n$ we write $|\vx|$ for the euclidean length of $\vx$.
\begin{definition}
Let $\M$ and $\Da>1$ be positive integers and let $L$ be a non-negative real.
We say that a set $S$ is in Lip$(\Da,\M,L)$ if 
$S$ is a subset of $\IR^\Da$, and 
if there are $\M$ maps 
$\phi_1,...,\phi_M:[0,1]^{\Da-1}\longrightarrow \IR^\Da$
satisfying a Lipschitz condition
\begin{alignat}3
\label{lipcond1}
|\phi_i(\vx)-\phi_i(\vy)|\leq L|\vx-\vy| \text{ for } \vx,\vy \in [0,1]^{\Da-1}, i=1,...,M 
\end{alignat}
such that $S$ is covered by the images
of the maps $\phi_i$.
\end{definition}
We call $L$ a Lipschitz constant for the maps $\phi_i$. By definition the empty set
lies in Lip$(\Da,\M,L)$ for any positive integers $\M$ and $\Da>1$ and any 
non-negative $L$.
\begin{definition}[Adelic-Lipschitz system]\label{defALS}
An adelic-Lipschitz system ($\ALS$) 
$\en_K$ or simply $\en$ on $K$ (of dimension $n$) is
a set of continuous maps
\begin{alignat}3
\label{Abb1}
N_v: K_v^{n+1}\rightarrow \Gamma_v \quad v \in M_K
\end{alignat}
such that for $v \in M_K$ we have
\begin{alignat*}3
(i)&\text{ } N_v({\vz })=0 \text{ if and only if } {\vz} ={\vnull},\\
(ii)&\text{ } N_v(\omega {\vz})=|\omega|_v N_v({\vz}) \text{ for all
$\omega$ in $K_v$ and all ${\vz}$ in $K_v^{n+1}$},\\
(iii)&\text{ if $v \mid \infty$: }\{{\vz};N_v({\vz})=1\} \text{
is in $Lip(d_v(n+1),\M_v,L_v)$ for some $\M_v, L_v$},\\
(iv)&\text{ if $v \nmid \infty$: }N_v({\vz_1}+{\vz_2})
\leq \max\{N_v({\vz}_1),N_v({\vz}_2)\} \text{ for all 
${\vz}_1,{\vz}_2$ 
in $K_v^{n+1}$}.
\end{alignat*}
Moreover we assume that 
\begin{alignat}3
\label{Nvmaxnorm}
N_v(\vz)=\max\{|z_0|_v,...,|z_n|_v\}
\end{alignat}
\end{definition}
for all but a finite number
of $v \in M_K$.
If we consider only the functions $N_v$ for $v\mid\infty$
then we get an $(r,s)$-Lipschitz system (of dimension $n$)
in the sense of Masser and Vaaler \cite{1}.
With $\M_v$ and $L_v$ from $(iii)$ we define
\begin{alignat*}1
\M_{\en}&=\max_{v\mid \infty}\M_v,\\
L_{\en}&=\max_{v\mid \infty}L_v.
\end{alignat*}
We say that $\en$
is an $\ALS$ with associated constants $\M_{\en}, L_{\en}$.
The set defined in $(iii)$ is the boundary 
of the set ${\bf B}_v=\{{\vz};N_v({\vz})<1\}$
and therefore ${\bf B}_v$ is a bounded symmetric open star-body
in $\IR^{n+1}$ or $\IC^{n+1}$ (see also \cite{1} p.431). In particular ${\bf B}_v$ has a finite volume $V_v$.\\

Let us consider the system
where $N_v$ is as in (\ref{Nvmaxnorm}) for all places $v$.
If $v$ is an infinite place then  
${\bf B}_v$ is a
cube for $d_v=1$ and the complex analogue
if $d_v=2$. Their boundaries are clearly
in Lip$(d_v(n+1),M_v,L_v)$ most naturally
with $M_v=2n+2$ maps and $L_v=2$
if $d_v=1$ and 
with $M_v=n+1$ maps and for example $L_v=2\pi\sqrt{2n+1}$
if $d_v=2$.
This system is called the standard
adelic-Lipschitz system.\\

We return to arbitrary adelic-Lipschitz systems.
We claim that for any $v\in M_K$ there is a $c_v$ in the value group
$\Gamma_v^*=\Gamma_v\backslash\{0\}$ with
\begin{alignat}3
\label{Nineq1}
N_v({\vz})\geq c_v\max\{|z_0|_v,...,|z_n|_v\}
\end{alignat}
for all $\vz=(z_0,...,z_n)$ in $K_v^{n+1}$.
For if  $v$ is archimedean then ${\bf B}_v$ is
bounded open and it contains the origin.
Since $\Gamma_v^*$ contains arbitrary small
positive numbers the
claim follows by $(ii)$.
Now for $v$ non-archimedean $N_v$ and $\max\{|z_0|_v,...,|z_n|_v\}$ define norms on 
the vector space $K_v^{n+1}$ over the complete field $K_v$.
But on a finite dimensional vector space over a complete
field all norms are equivalent (\cite{2} Corollary 5. p.93)
hence (\ref{Nineq1}) remains true for a suitable choice
of $c_v$.\\

So let $\en$ be an $\ALS$ on $K$ of dimension $n$. For
every $v$ in $M_K$ let
$c_v$ be an element of $\Gamma_v^*$,
such that $c_v\leq 1$ and (\ref{Nineq1}) holds.
Due to (\ref{Nvmaxnorm}) we can assume that $c_v\neq 1$ only for a finite number of places $v$.
We define
\begin{alignat}3
\label{defcfin}
C^{fin}_{\en}&=\prod_{v\nmid \infty}c_v^{-\frac{d_v}{d}}\geq 1
\end{alignat}
and
\begin{alignat}3
\label{defcinf}
C^{inf}_{\en}&=\max_{v\mid \infty}\{c_v^{-1}\}\geq 1.
\end{alignat}
Multiplying the finite and the infinite part 
gives rise to another constant 
\begin{alignat}3
\label{defc}
C_{\en}&=C^{fin}_{\en}C^{inf}_{\en}.
\end{alignat}
It will turn out that besides $\M_{\en}$ and $L_{\en}$ this is another important quantity for an $\ALS$. So we say that \it $\en$
is an $\ALS$ with associated constants $C_{\en},\M_{\en},L_{\en}$.\rm 
\begin{remark}\label{normconvex}
Let $v$ be an infinite place.
Suppose $N_v:K_v^{n+1}\longrightarrow [0,\infty)$
defines a norm, so that 
$N_v({\vz_1}+{\vz_2})
\leq N_v({\vz}_1)+N_v({\vz}_2)$. Then
${\bf B}_v$ is convex and (\ref{Nineq1}) combined with 
(\ref{defcfin}), (\ref{defcinf}) and (\ref{defc}) shows that ${\bf B}_v$ lies in
$B_0(C_{\en}\sqrt{n+1})$.
This implies (see Theorem A.1 in \cite{WiThesis}) that  $\partial{\bf B}_v$ lies in Lip$(d_v(n+1),1,8{d_v}^2(n+1)^{5/2}C_{\en})$.
\end{remark}

We denote by $\sigma_1,...,\sigma_d$ 
the embeddings from $K$ to $\IR$ or
$\IC$ respectively, ordered such that
$\sigma_{r+s+i}=\overline{\sigma}_{r+i}$ for 
$1\leq i \leq s$.
We define
\begin{alignat}3
\label{sigd}
&\sigma:K\longrightarrow \IR^r\times\IC^s\\
\nonumber&\sigma(\alpha)=(\sigma_1(\alpha),...,\sigma_{r+s}(\alpha)).
\end{alignat}
Sometimes it will be more readable to omit the brackets
and simply to write $\sigma\alpha$.
We identify $\IC$ in the usual way with $\IR^2$
and extend $\sigma$ componentwise to get a map
\begin{alignat}3
\label{sigD}
\sigma:K^{n+1}\longrightarrow \IR^{D}
\end{alignat}
where $D=d(n+1)$.
On $\IR^{D}$ we use $|\cdot|$ for the
usual euclidean norm. For $v\in M_K$
let $\sigma_v$ be the canonical embedding of $K$ in $K_v$,
again extended componentwise on $K^{n+1}$. 
\begin{definition}
Let $\D\neq 0$ be a fractional ideal in $K$ and $\en$
an $\ALS$ of dimension $n$. We define
\begin{alignat}3
\label{defLamen}
\Lamen(\D)=\{\sigma(\balf); \balf \in K^{n+1}, 
N_v(\sigma_v\balf)\leq |\D|_v \text{ for all finite }v \}
\end{alignat}
where $|\D|_v=N\pw_v^{-\frac{\ord_{\pw_v}\D}{d_v}}$.
\end{definition}
It is easy to see that 
$\Lamen(\D)$ is an additive subgroup of $\IR^D$.
Now assume $B\geq 1$ and 
$|\sigma(\balf)|\leq B$; then (\ref{Nineq1}) implies
$H(\balf)^d\leq (BC_{\en}^{fin})^d N\D^{-1}$ and by Northcott's Theorem
we deduce that $\Lamen(\D)$ is discrete.
The same argument as for (\ref{Nineq1}) yields positive real numbers $C_v$, one for each non-archimedean place $v\in M_K$,
with $N_v({\vz})\leq C_v\max\{|z_0|_v,...,|z_n|_v\}$ for all $\vz=(z_0,...,z_n)$ in $K_v^{n+1}$
and $C_v=1$ for all but finitely many non-archimedean $v\in M_K$.
Thus there exists an ideal $\C_1\neq 0$ in $\Oseen_K$ with 
$|\C_1|_v\leq 1/C_v$ for all non-archimedean places $v\in M_K$.
This means that $\sigma(\C_1\D)^{n+1}\subseteq \Lamen(\D)$.
It is well-known that the additive group $\sigma(\C_1\D)^{n+1}$ has maximal rank in $\IR^D$.
Therefore $\Lamen(\D)$ is  a discrete additive subgroup of $\IR^D$
of maximal rank. Hence $\Lamen(\D)$
is a lattice.
Notice that for $\varepsilon$ in $K^*$ one has
\begin{alignat}3
\label{deltawelldef1}
\det\Lamen((\varepsilon)\D)=
|N_{K/\IQ}(\varepsilon)|^{n+1}\det\Lamen(\D).
\end{alignat}
Therefore
\begin{alignat}3
\label{idkl}
\Delta_{\en}(\mathcal{D})=\frac{\det\Lambda_{\en}(\D)}
{N\D^{n+1}}
\end{alignat}
is independent of the choice of the representative $\D$
but depends only on the  ideal class $\mathcal{D}$ of $\D$. 
Let $\Cl_K$ denote the ideal class group of $K$.
We define
\begin{alignat}3
\label{defVfin1}
V_{\en}^{fin}=2^{-s(n+1)}|\Delta_K|^{\frac{n+1}{2}}h_K^{-1}
\sum_{\mathcal{D} \in \Cl_K}\Delta_{\en}(\mathcal{D})^{-1}
\end{alignat}
for the finite part, where as usual,
$s$ denotes the number of pairs of complex conjugate
embeddings of $K$, $h_K$ the class number 
of $K$ and $\Delta_K$ the discriminant of $K$.
The infinite part is defined by 
\begin{alignat*}3
V_{\en}^{inf}=\prod_{v \mid \infty}V_{v}.
\end{alignat*}
By virtue of (\ref{Nineq1}) we observe that
\begin{alignat}3
\label{Vinfbou1}
V_{\en}^{inf}=\prod_{v|\infty} V_v\leq 
\prod_{v|\infty}(2 C^{inf}_{\en})^{d_v(n+1)}=
(2 C^{inf}_{\en})^{d(n+1)}.
\end{alignat} 
We multiply the finite and the infinite part
to get a global volume 
\begin{alignat}3
\label{defVen}
V_{\en}=V_{\en}^{inf}V_{\en}^{fin}.
\end{alignat}
\subsection{Adelic-Lipschitz heights on $\mathbb{P}^n(K)$}\label{2subsec2}
Let $\en$ be an $\ALS$ on $K$
of dimension $n$. Then the height $\hen$ 
on $K^{n+1}$ is defined by
\begin{alignat*}3
\hen(\balf)=\prod_{v \in M_K} N_v(\sigma_v(\balf))^{\frac{d_v}{d}}.
\end{alignat*}
Thanks to the product formula
and $(ii)$ from Subsection \ref{subsecdefALS} $\hen(\balf)$ does not change
if we multiply each coordinate of $\balf$ with a fixed element of $K^*$.
Therefore $\hen$ is well-defined on $\IP^n(K)$ by setting
\begin{alignat*}3
\hen(P)=\hen(\balf)
\end{alignat*}
where $P=(\alpha_0:...:\alpha_n) \in \IP^n(K)$ and $\balf=(\alpha_0,...,\alpha_n) \in K^{n+1}$.
\begin{remark}
Multiplying (\ref{Nineq1}) over all places with 
suitable multiplicities yields 
\begin{alignat}3
\label{HAquiv}
\hen(P)\geq C_{\en}^{-1} H(P)
\end{alignat}
for $P\in \IP^n(K)$.
Thanks to Northcott's Theorem it follows that
$\{P\in \IP^n(K); \hen(P)\leq X\}$ is a finite set for
each $X$ in $[0,\infty)$.
\end{remark}

\subsection{Adelic-Lipschitz systems on a collection of number fields}\label{subsec1.1}
Recall that $k$ is a number field of degree $\m$ and $\kbar$ is an algebraic closure of $k$. 
We fix $k$ and $\kbar$ throughout and assume finite extensions of $k$ to lie in $\kbar$.
Let $\coll$ be a collection of finite extensions of $k$. We
are especially interested in the set of all extensions
of fixed relative degree. We denote it by
\begin{alignat*}1
\coll_e=\coll_e(k)=\{K\subseteq \kbar;[K:k]=e\}.
\end{alignat*}
Let 
$\en$ be a collection of adelic-Lipschitz systems $\en_K$ of dimension $n$ - one for each $K$ of $\coll$. Then we call $\en$ an \em adelic-Lipschitz system $(\ALS)$ on 
$\coll$ of dimension $n$. \rm
We say $\en$ is a \em uniform \rm $\ALS$ on $\coll$
of dimension $n$
with associated constants $C_{\en},\M_{\en},L_{\en}$ in $\IR$ if
the following holds: 
for each $\ALS$ $\en_K$ of the collection $\en$
we can choose associated constants
$C_{\en_K},\M_{\en_K},L_{\en_K}$ satisfying
\begin{alignat*}1
C_{\en_K}\leq C_{\en},\quad
\M_{\en_K}\leq \M_{\en},\quad
L_{\en_K}\leq L_{\en}.
\end{alignat*}
Notice that a uniform $\ALS$ $\en$ (of dimension $n$) on the collection consisting only of a single field $K$ with associated constants $C_{\en},\M_{\en},L_{\en}$ is simply an $\ALS$
$\en$ (of dimension $n$) on $K$ with associated constants
$C_{\en},\M_{\en},L_{\en}$ in the sense of Subsection \ref{subsecdefALS}.\\

A standard example for a uniform $\ALS$ on $\coll_e$ (of dimension $n$) is given as 
follows: for each $K$ in $\coll_e$ choose the standard $\ALS$ on $K$ (of dimension $n$)
so that $N_v$ is as in (\ref{Nvmaxnorm}) for each
$v$ in $M_K$. For this system we may choose
$C_{\en}=1$,
$\M_{\en}=2n+2$ and $L_{\en}=2\pi\sqrt{2n+1}$.
Choosing $l^2$-norms at all infinite places
and $N_v$ as in (\ref{Nvmaxnorm}) for all finite places yields another important uniform  $\ALS$.

\subsection{Adelic-Lipschitz heights on $\mathbb{P}^n(k;e)$}\label{ALHoncoll}
Let $\coll$ be a collection of finite extensions of $k$ and let $\en$ be an $\ALS$ of dimension $n$ on $\coll$.
Now we can define heights on $\IP^n(K/k)$ (the set of points $P\in\IP^n(K)$ with $k(P)=K$)
for $K$ in $\coll$.
Let $P\in \IP^n(K/k)$.
According to Subsection \ref{2subsec2} we know that $H_{\en_{K}}(\cdot)$
defines a projective height on $\IP^n(K)$.
Now we define
\begin{alignat}1
\label{defALH2}
\hen(P)=H_{\en_{K}}(P).
\end{alignat}
From Subsection \ref{2subsec2} we know
\begin{alignat}1
\label{defALH3}
H_{\en_{K}}(P)=\prod_{v\in M_{K}} N_v(\sigma_v(\balf))^{\frac{d_v}{d}}
\end{alignat}
for the functions $N_v$ of $\en_{K}$ and
$[K:\IQ]=d$, $[K_v:\IQ_v]=d_v$. Starting with $\coll=\coll_e$
we get a height defined on $\IP^n(k;e)$.

\section{The main result}\label{introchap4}
Let $\en$ be an $\ALS$ on $\coll_e$ of dimension $n$.
Then $\hen(\cdot)$ defines a height on $\IP^n(k;e)$, the set of points 
$P=(\alpha_0:...:\alpha_n)$ in $\IP^n(\kbar)$ with $[k(P):k]=e$
where $k(P)=k(...,\alpha_i/\alpha_j,...)$
for $0\leq i,j\leq n$, $\alpha_j\neq 0$.
The associated counting function $Z_{\en}(\IP^n(k;e),X)$ denotes the number of points 
$P$ in $\IP^n(k;e)$ with $\hen(P)\leq X$. Assume $\en$ is a
uniform $\ALS$ on $\coll_e$
(of dimension $n$). Then due to Northcott and 
(\ref{HAquiv}) $Z_{\en}(\IP^n(k;e),X)$
is finite for all $X$ in $[0,\infty)$. 
The Schanuel constant $S_K(n)$ is defined as follows
\begin{alignat}1\label{Schanuelconst}
S_K(n)=\frac{h_KR_K}{w_K\zeta_K(n+1)}
(\frac{2^{r_K}(2\pi)^{s_K}}{\sqrt{|\Delta_K|}})^{n+1}
(n+1)^{r_K+s_K-1}.
\end{alignat}
Here $h_K$ is the class number, $R_K$ the regulator,
$w_K$ the number of roots of unity in $K$, $\zeta_K$
the Dedekind zeta-function of $K$, $\Delta_K$ the discriminant,
$r_K$ is the number of real embeddings of $K$ and $s_K$ is
the number of pairs of distinct complex conjugate embeddings of $K$.
Recall also the definition of $V_{\en_K}$ (see (\ref{defVen})).
Now we define the sum
\begin{alignat}1
\label{kkonst}
\Ce_{\en}=\Ce_{\en}(k,e,n)=\sum_{K}2^{-r_K(n+1)}\pi^{-s_K(n+1)}V_{\en_K}S_K(n)
\end{alignat}
where the sum runs over all extensions of $k$ 
with relative degree $e$.
We will prove that the sum in (\ref{kkonst}) 
converges if $n$ is large enough compared to $e$.
It will be convenient to use Landau's $O$-notation.
For non-negative real functions $f(X), g(X), h(X)$ we say that
$f(X)=g(X)+O(h(X))$ as $X>X_0$ tends to infinity if there is a constant $C_0$ such that
$|f(X)-g(X)|\leq C_0h(X)$ for each $X>X_0$.\\

After all this we are ready to state the main result.

\begin{theorem}\label{main theorem}
Let $e,n$ be positive integers and $k$ a number field
of degree $\m$. Suppose $\en$ is a uniform adelic-Lipschitz system of dimension $n$
on $\coll_e$, the collection of all finite extensions of $k$ of relative degree $e$, with associated constants $C_{\en},\M_{\en}$ and $L_{\en}$. Write
\begin{alignat}1\label{defAN}
A_{\en}&=\M_{\en}^{\m e}(C_{\en}(L_{\en}+1))^{\m e(n+1)-1}.
\end{alignat}
Suppose that either $e=1$ or
\begin{alignat}1\label{necond} 
n>{5e}/{2}+4+2/(\m e).
\end{alignat}
Then the sum in (\ref{kkonst}) converges and 
as $X>0$ tends to infinity we have
\begin{alignat*}1
Z_{\en}(\IP^n(k;e),X)=\Ce_{\en}X^{\m e(n+1)}
+O(A_{\en}X^{\m e(n+1)-1}\L),
\end{alignat*}
where $\L=\log\max\{2,2C_{\en}X\}$ if
$(\m e,n)=(1,1)$ and $\L=1$ otherwise.
The constant in $O$ depends only on
$k,e$ and $n$.
\end{theorem}

In subsequent papers \cite{art3} and \cite{art4} we will explore some
applications of Theorem \ref{main theorem}. Here we are content
with some immediate consequences.
For $e=1$ we recover a version of the
Proposition in \cite{1}, which allows more general norms at the finite places (this generalization will
be essential to deduce the results of \cite{art3}).
Now choose the standard uniform $\ALS$ as described at the end of Subsection \ref{subsec1.1} so that $\hen$ becomes the Weil height. Schanuel's Theorem
implies $S_K(n)=\Ce_{\en}(K,1,n)=2^{-r_K(n+1)}\pi^{-s_K(n+1)}V_{\en_K}S_K(n)$. We can verify 
\begin{alignat}1
\label{VenWeilheight}
V_{\en_K}=2^{r_K(n+1)}\pi^{s_K(n+1)}
\end{alignat} 
directly by noting
that $\Lamen(\D)=(\sigma\D)^{n+1}$ in (\ref{defLamen}),
so that $\det \Lamen(\D)=(2^{-s_K}N\D\sqrt{|\Delta_K|})^{n+1}$
(see \cite{1} Lemma 5).
Inserting the latter in definition (\ref{defVfin1}) we get $V_{\en_K}^{fin}=1$ and it is clear that
$V_{\en_K}^{inf}=2^{r_K(n+1)}\pi^{s_K(n+1)}$.
Then (\ref{VenWeilheight}) follows from $V_{\en_K}=V_{\en_K}^{inf}V_{\en_K}^{fin}$
and so we find Theorem \ref{mainthintro} from the introduction.
For $k=\IQ$ and $e=2$ we recover essentially  Schmidt's theorem (\ref{ThSchm2}) but only for $n>10$ while Schmidt does it for all
$n\geq 3$ and even (in a modified form) for $n=1,2$. For $k=\IQ$ and $e>2$ we find Gao's result (\ref{ThGao})
but again with the stronger restriction $n>5e/2+4+2/(\m e)$
instead of Gao's $n>e$.\\

It is likely that Theorem \ref{main theorem} is valid for $n>e$ instead of (\ref{necond}).
Gao showed, at least for his definition of height (see also 
\cite{WiThesis} Appendix B),
that for $k=\IQ$
the bound $n>e$ suffices. On the other hand Schmidt's
lower bound in (\ref{ThSchm1'}) implies that
Theorem \ref{main theorem} cannot hold for $e>1$ and $n<e$.
However, there is a good possibility of obtaining
the asymptotics for $e>1$ and $n=1$ using a kind
of generalized Mahler measure.

\section{Proof of the main result}\label{2subsec5}
The major part of the work was already done in \cite{art1} where 
we proved estimates for $Z_{\en_K}(\IP^n(K/k),X)$. These estimates will be essential to deduce 
Theorem \ref{main theorem}.
\subsection{Preliminaries}\label{prelim}
Let $K$ be in $\coll_e$. Then by definition $\hen(P)=\henK(P)$
for all $P$ in $\IP^n(K/k)$.
Since 
\begin{alignat}1
\label{Gl1.5.1}
\IP^n(k;e)=\bigcup_{K\in \coll_e}\IP^n(K/k)
\end{alignat}
where the right hand side is a disjoint union,
we get 
\begin{alignat}1
\label{Gl1.5.2}
Z_{\en}(\IP^n(k;e),X)=
\sum_{K \in \coll_e}Z_{\en_K}(\IP^n(K/k),X).
\end{alignat}
To state the estimates for $Z_{\en_K}(\IP^n(K/k),X)$ from \cite{art1} we are forced to introduce some more notation.
For fields $k,K$ with $k\subseteq K$ and $[K:k]=e$ we define
\begin{alignat*}3
G(K/k)=
\{[K_0:k]; \text{$K_0$ is a field with $k\subseteq K_0\subsetneq K$}\}
\end{alignat*}
if $k\neq K$, and we define  
\begin{alignat*}3
G(K/k)=\{1\}
\end{alignat*}
if $k=K$. Clearly $|G(K/k)|\leq e$.
Then for an integer $g\in G(K/k)$ we define
\begin{alignat}3
\label{delta2}
\delta_g(K/k)=
\underset{\alpha,\beta}\inf\{H(1,\alpha,\beta); k(\alpha,\beta)=K,
[k(\alpha):k]=g\}
\end{alignat}
(which is $\geq 1$) and 
\begin{alignat}3
\label{mu2}
\mu_g=m(e-g)(n+1)-1.
\end{alignat}
In \cite{art1} the author proved the following result.
\begin{theorem}\label{prop3}
Let $k,K$ be number fields with $k\subseteq K$ and $[K:k]=e$, $[k:\IQ]=m$, $[K:\IQ]=d$. Let $\en_K$ be an adelic-Lipschitz system of dimension $n$ 
on $K$ with associated constants $C_{\en_K},L_{\en_K},\M_{\en_K}$.
Write
\begin{alignat}3
\label{defANK}
A_{\en_K}&=\M_{\en_K}^{d}(C_{\en_K}(L_{\en_K}+1))^{d(n+1)-1},\\
\label{mainterm1}
\Mainterm&=2^{-r_K(n+1)}\pi^{-s_K(n+1)}V_{\en_K}S_K(n),\\
\label{errorterm1}
\Errorterm&=A_{\en_K}R_Kh_K\sum_{g\in G(K/k)}\delta_{g}(K/k)^{-\mu_g}.
\end{alignat}
Then as $X>0$ tends to infinity we have
\begin{alignat*}3
Z_{\en_K}(\IP^n(K/k),X)=
\Mainterm X^{d(n+1)}
+O(\Errorterm X^{d(n+1)-1}\L),
\end{alignat*}
where
\begin{alignat*}3
\L&=\log\max\{2,2C_{\en_K}X\} \text{ if }(n,d)=(1,1)\text{ and }\L=1 \text{ otherwise}
\end{alignat*}
and the implied constant in $O$ depends only on $n$ and $d$.
\end{theorem}
Thanks to (\ref{Gl1.5.2}) and Theorem \ref{prop3} it 
suffices to show that $\sum \Mainterm$ and $\sum \Errorterm$
are convergent
(here the sum runs over the same fields as in (\ref{kkonst}) and
(\ref{Gl1.5.2})).\\

We will also deal with $\delta(\cdot)$, a simplified version of $\delta_g(\cdot)$ 
\begin{alignat*}3
\delta(K/k)=\underset{\alpha}\inf\{H(1,\alpha);K=k(\alpha)\}.
\end{alignat*}
The quantity $\delta(K/\IQ)$ was already introduced by Roy and Thunder \cite{8}.\\

It will be extremely convenient to use Vinogradov's $\ll,\gg$-notation.
The constants involved in $\ll$ and $\gg$ will depend only on
$k,n,e$ unless we indicate the dependence on additional parameters by an index.\\

The case $e=1$ of Theorem \ref{main theorem} is already covered by Theorem \ref{prop3} by choosing $K=k$.
For the rest of this article we assume 
\begin{alignat*}1
e>1.
\end{alignat*}
For a non-zero ideal $\A$ in $K$ let $D_{K/k}(\A)$
be the discriminant-ideal of $\A$ relative to $k$ (for definitions see \cite{23} p.212 or \cite{13})
and write $D_{K/k}$ for $D_{K/k}(\Oseen_K)$  where $\Oseen_K$ denotes the ring of integers in $K$.
By assumption we have $\IQ\subseteq k \subseteq K$ and hence by \cite{23} (2.10) Korollar p.213
\begin{alignat}1
\label{dis}
|\Delta_{K/\IQ}|=|\Delta_{k/\IQ}|^{[K:k]}N_{k/\IQ}(D_{K/k})
\end{alignat}
where $N_{k/\IQ}(a)$ denotes the absolute norm of an ideal $a\neq 0$ of the ring of integers
$\Oseen_k$, i.e. $N_{k/\IQ}(a)=|\Oseen_k/a|$.
Let $P$ be in $\IP^n(K/k)$, so
$K=k(P)$. We use a theorem of Silverman (\cite{9} Theorem 2) 
with Silverman's $S_F$ (for $F=k$) as the set of archimedean absolute values.
Then Silverman's $L_F(\cdot)$ is simply the usual norm
$N_{k/\IQ}(\cdot)$. Hence we deduce
\begin{alignat}1
\label{sil}
H(P)^{\m }\geq \exp \left(-\frac{\delta_k\log e}{2(e-1)}\right)
N_{k/\IQ}(D_{K/k})^{\frac{1}{2e(e-1)}}
\end{alignat}
where $\delta_k$ is the number of archimedean places
in $M_k$. Since Silverman uses not an absolute height
but rather an ``absolute height relative to $k$'', we had 
to take the $\m $-th power on the left hand side of 
(\ref{sil}).
Combining (\ref{dis}) and (\ref{sil}) yields
\begin{alignat}1
\label{sil2}
H(P)&\geq \exp \left(-\frac{\delta_k\log e}{2(e-1)\m}\right)
|\Delta_k|^{-\frac{1}{2(e-1)\m }}
|\Delta_K|^{\frac{1}{2e(e-1)\m }}\\
\nonumber&\gg |\Delta_K|^{\frac{1}{2e(e-1)\m }}.
\end{alignat}
Recalling the definitions of $\delta$, $\delta_g$ and
$G(K/k)$ and taking $P=(1:\alpha_1:\alpha_2)$ in $\IP^2(K/k)$ we get 
\begin{alignat}1
\label{de1}
\delta_{g}(K/k)
\gg |\Delta_K|^{\frac{1}{2e(e-1)\m }}
\end{alignat}
for any $g\in G(K/k)$; and similarly
\begin{alignat}1
\label{de11}
\delta(K/k)
\gg |\Delta_K|^{\frac{1}{2e(e-1)\m }}.
\end{alignat}
Here it might be worthwile to point out that 
(\ref{sil2}) can be used to prove 
a version of Theorem \ref{prop3} where 
$\Errorterm$ is redefined in terms of the discriminants; namely
\begin{alignat}1
\label{discerrorterm}
\Errorterm=A_{\en_K}R_Kh_K\sum_{g\in G(K/k)}(|\Delta_k|^{-e}|\Delta_K|)^{-\frac{\mu_g}{2e(e-1)\m}}.
\end{alignat}
At a first glance this error term looks more appropriate
due to the unavoidable appearance of $\Delta_K$ in the main term.
But as it turns out, the summation over $\Delta_K$ instead of over $\delta_g(K/k)$ leads to a result
weaker than Theorem \ref{main theorem}, in which we have to assume that $n$ exceeds some quadratic function of $e$ instead
of (\ref{necond}). The reason for this is, that we have estimates for the number of number fields $K$ with $\delta_g(K/k)\leq T$
which are more accurate than the best available estimates for the number of number fields with $|\Delta_K|\leq T$, see Section \ref{countingfields}
for a discussion on this.
Thanks to the well-known Theorem of Siegel-Brauer
(\cite{13} p.328 Corollary or \cite{24} p.67 Satz 1 for 
a more precise version) we
can use the inequalities (\ref{de1}) and (\ref{de11}) to bound the product of regulator and class number.
More precisely we have
\begin{alignat}1
\label{de2}
R_Kh_K \ll_{\beta} \delta_{g}(K/k)^{\beta}
\end{alignat}
and 
\begin{alignat}1
\label{de22}
R_Kh_K \ll_{\beta} \delta(K/k)^{\beta}.
\end{alignat}
for any $\beta>e(e-1)\m $ and any $g\in G(K/k)$.

\subsection{Three preparatory lemmas}
We start with a very simple argument, known as
dyadic summation. Since it will be frequently used we state it as a lemma.
\begin{lemma}[Dyadic summation]\label{dyadicsum}
Let $\coll$ be a non-empty subset of $\coll_e$ and let $\f$ be a map $\f:\coll \longrightarrow [1,\infty)$. 
Write $N_{\f}(T)=|\{K\in \coll;\f(K)\leq T\}|$ and
suppose there
are nonnegative real numbers $b,c$ (independent of $T$) 
with 
\begin{alignat*}1
N_{\f}(T)\leq cT^b
\end{alignat*}
for every $T>0$.
Let $\coll'$ be a non-empty subset of $\coll$.
Set $\mathfrak{M}=[\log_2\max_{\coll'}\f(K)]+1$ if 
$\coll'$ is finite and $\mathfrak{M}=\infty$ 
otherwise.
Moreover
suppose $\alpha$ is a real number such that
$\sum_{i=1}^{\mathfrak{M}}2^{i(\alpha+b)}$ converges.
Then we have
\begin{alignat*}1
\sum_{K\in \coll'}\f(K)^{\alpha}\leq c 2^{|\alpha|}\sum_{i=1}^{\mathfrak{M}}2^{i(\alpha+b)}.
\end{alignat*}
\end{lemma}
\begin{rproof}
From the definition of $\mathfrak{M}$ and since $\coll'\subseteq \coll$ we have
\begin{alignat*}1
\sum_{K\in \coll'}\f(K)^{\alpha}=
\sum_{i=1}^{\mathfrak{M}}\sum_{K\in \coll' \atop 2^{i-1}\leq \f(K)<2^i}
\f(K)^{\alpha}\leq \sum_{i=1}^{\mathfrak{M}}\sum_{K\in \coll \atop 2^{i-1}\leq \f(K)<2^i}
\f(K)^{\alpha}.
\end{alignat*}
First suppose $\alpha<0$. Then the latter is
\begin{alignat*}1
\leq \sum_{i=1}^{\mathfrak{M}}2^{(i-1)\alpha}N_{\f}(2^i)
\leq
c 2^{-\alpha}\sum_{i=1}^{\mathfrak{M}}2^{i(\alpha+b)}.
\end{alignat*}
If $\alpha\geq 0$ then we even get
\begin{alignat*}1
\sum_{K\in \coll'}\f(K)^{\alpha}\leq
c\sum_{i=1}^{\mathfrak{M}}2^{i(\alpha+b)}.
\end{alignat*}
This proves the lemma.
\end{rproof}
Recall the definition of $G(K/k)$ from Subsection \ref{prelim}.
In our applications $\f$ will be $\delta_g$ and
$\coll$ will be         
\begin{alignat*}1
\collg=\{K\in \coll_e; g\in G(K/k)\}
\end{alignat*}
the set of extensions 
$K$ of $k$ of relative degree $e$
containing an intermediate field $K_0\subsetneq K$ with $[K_0:k]=g$.
Let $G_u$ be the union of all sets $G(K/k)$
\begin{alignat*}1
G_u=\bigcup_{K\in \coll_e}G(K/k),
\end{alignat*}
so that $\collg$ is non-empty if and only if $g \in G_u$.
In fact $G_u$ is simply the set of positive, proper
divisors of $e$ but we need only 
\begin{alignat*}1
\{1\}\subseteq G_u \subseteq \{1,...,[e/2]\}.
\end{alignat*}
To apply the dyadic summation lemma we need information
about the growth rate of $N_{\delta_g}(T)$.
In accordance with the notation in Lemma \ref{dyadicsum} we define for an integer $g\in G_u$ 
and real positive $T$ 
\begin{alignat*}1
N_{\delta_g}(T)=|\{K\in \collg;\delta_g(K/k)\leq T\}|.
\end{alignat*}
The set on the right-hand side is finite. More precisely
we have the following lemma.
\begin{lemma}\label{lemmaNdeltas}
Set $\gamma_g=\m (g^2+g+e^2/g+e)$. Then for real positive
$T$ and $g$ in $G_u$ we have
\begin{alignat*}1
N_{\delta_g}(T)\ll T^{\gamma_g}.
\end{alignat*}
\end{lemma}
\begin{rproof}
Since $H(1,\alpha_1,\alpha_2)\geq \max\{H(1,\alpha_1),H(1,\alpha_2)\}$ it suffices to show that the number of tuples $(\alpha_1,\alpha_2)\in \kbar^2$ with
\begin{alignat}1
\label{conddeg1}
&[k(\alpha_1):k]=g\\
\label{conddeg2}
&[k(\alpha_1,\alpha_2):k(\alpha_1)]=e/g\\
\label{condheight}
&H(1,\alpha_1), H(1,\alpha_2)\leq T
\end{alignat}
is $\ll T^{\gamma_g}$.
The number of projective points
in $\IP(k;g)$ with height not exceeding $T$ is an upper
bound for the number of $\alpha_1$ in $\kbar$ of relative degree
$g$ with $H(1,\alpha_1)\leq T$. 
Thus by  (\ref{ThSchm1'}) we get the upper bound 
\begin{alignat}1
\label{ubalpha1}
\ll T^{\m g(g+1)}
\end{alignat}
for the number of $\alpha_1$.
Next for each $\alpha_1$ we count the number of $\alpha_2$.
By (\ref{conddeg2}) we have 
$[k(\alpha_1,\alpha_2):k(\alpha_1)]=e/g$
and moreover $H(1,\alpha_2)\leq T$. Applying (\ref{ThSchm1'})
(note that the constant
$\CS(k,e,n)$ in (\ref{ThSchm1'}) depends only on $[k:\IQ],e,n$)
once more yields the upper bound
\begin{alignat}1
\label{ubalpha2}
\ll T^{[k(\alpha_1):\IQ](e/g)(e/g+1)}=T^{\m e(e/g+1)}
\end{alignat}
for the number of $\alpha_2$ provided $\alpha_1$ is fixed.
Multiplying the bound (\ref{ubalpha1}) for the number
of $\alpha_1$ and (\ref{ubalpha2}) gives the upper bound
\begin{alignat*}1
\ll T^{\m(g^2+g+e^2/g+e)}
\end{alignat*}
for the number of tuples $(\alpha_1,\alpha_2)$ and thereby
proves the lemma.
\end{rproof}
Recall that $\delta_1=\delta$
and that $N_{\delta}(T)$ denotes the number of number fields $K$
in $\kbar$ of relative degree $e$ with $\delta(K/k)\leq T$.
So Lemma \ref{lemmaNdeltas} with $g=1$ yields an upper bound for the growth rate of
$N_{\delta}(T)$ but applying (\ref{ThSchm1'}) directly gives a slightly better result.
\begin{lemma}\label{lemmagamma}
Set $\gamma=\m e(e+1)$ and let $C_{\delta}=\CS(k,e,1)$ be as in 
(\ref{ThSchm1'}). Then for $T>0$ we have
\begin{alignat}1
\label{Ndelta}
N_{\delta}(T)\leq C_{\delta} T^{\gamma}.
\end{alignat}
\end{lemma}
\begin{rproof}
The number of points in $\IP(k;e)$ with height not
larger than $T$ is clearly an upper bound for
$N_{\delta}(T)$. Thus the lemma follows from
the upper bound in (\ref{ThSchm1'}).
\end{rproof}
In fact Lemma \ref{lemmaNdeltas}
would suffice to prove the full Theorem \ref{main theorem}, so
one could omit Lemma \ref{lemmagamma}. We did not
because $\gamma$ looks nicer than $\gamma_1$
and the proof above is essentially simply a reference.

\subsection{Proof of Theorem \ref{main theorem}}
Recall the definition of $\Errorterm$ and $\Mainterm$ from (\ref{errorterm1}) and (\ref{mainterm1}).
We have seen that it 
suffices to show that $\sum \Errorterm$ and  $\sum \Mainterm$
are convergent where the sums run over all fields in $\coll_e$.\\

Since $\en$ is a uniform $\ALS$ on $\coll_e$ with associated constants
$C_{\en},\M_{\en}$ and $L_{\en}$ we can assume that
\begin{alignat}1 
\label{CenCenk}
C_{\en_K}&\leq C_{\en},\\
\label{MenMenk}
\M_{\en_K}&\leq \M_{\en},\\
\label{LenLenk}
L_{\en_K}&\leq L_{\en}.
\end{alignat}
Hence by definition (\ref{defAN}) and (\ref{defANK})
\begin{alignat}1
\label{Aen}
A_{\en_K}\leq A_{\en}.
\end{alignat}
Let us now prove that $\sum_K \Errorterm$ converges.
We set $\beta=e(e-1)\m +1/8$. Using (\ref{de2})
and (\ref{Aen}) we get
\begin{alignat*}1
\sum_{K\in \coll_e} \Errorterm \ll A_{\en}\sum_{K\in \coll_e} \sum_{g\in G(K/k)}\delta_{g}(K/k)^{\beta-\mu_g}.
\end{alignat*}
Recall that $G_u=\bigcup_{\coll_e} G(K/k)$.
So the term on the right-hand side above is
\begin{alignat}1
\nonumber
&A_{\en}\sum_{g\in G_u}\sum_{K \in \coll_e \atop g\in G(K/k)}\delta_{g}(K/k)^{\beta-\mu_g}\\
\label{EKsum1}=&A_{\en}\sum_{g\in G_u}\sum_{K\in \collg}
\delta_{g}(K/k)^{\beta-\mu_g}
\end{alignat}
provided the sum converges. This will be verified in a 
moment (see (\ref{confirmconv})). 
Applying the dyadic summation lemma with 
$\f=\delta_g$ and $b=\gamma_g$ from 
Lemma \ref{lemmaNdeltas} we see that the latter is
\begin{alignat*}1
\ll A_{\en}\sum_{g\in G_u}\sum_{i=1}^{\infty}2^{i(\gamma_g+\beta-\mu_g)}.
\end{alignat*}
The next lemma will tell us that
$\gamma_{g}+\beta-\mu_g\leq -1/8$. Assuming this for
a moment we see that
the inner sum above
is $\ll 1$. Thus we derive the upper bound
\begin{alignat}1
\label{confirmconv}
\ll A_{\en}\sum_{g\in G_u}1
\ll A_{\en},
\end{alignat}
confirming that the whole sum in (\ref{EKsum1}) converges. This verifies the convergence
of $\sum_K \Errorterm$ under the hypothesis $\gamma_{g}+\beta-\mu_g\leq -1/8$ for all $g\in G_u$.
The following lemma shows that this hypothesis holds true.
Recall that we assume $e>1$ and therefore by our assumption
in Theorem \ref{main theorem} 
$n>5e/2+4+2/(\m e)$.
\begin{lemma}\label{gabemuisneg}
Let $g$ be in $G_u$. Then
\begin{alignat}1
\label{ineqgabemu}
\gamma_{g}+\beta-\mu_g\leq -\frac{1}{8}.
\end{alignat}
\end{lemma}
\begin{rproof}
Recall that $G_u\subseteq \{1,...,[e/2]\}$ and $\mu_g=\m (e-g)(n+1)-1$.
Write
\begin{alignat*}1
\IG(g)=\frac{1}{\m(e-g)}(\gamma_{g}+\beta+1).
\end{alignat*}
So (\ref{ineqgabemu}) claims that $m(e-g)(\IG(g)-(n+1))\leq -1/8$
for all $g\in G_u$.
Hence it suffices to show that
\begin{alignat*}1
\IG(g)-(n+1)\leq -\frac{1}{4\m e}
\end{alignat*}
for $1\leq g\leq e/2$. By definition 
\begin{alignat*}1
\IG(g)=\frac{g^2+g+e^2/g+e}{e-g}+
\frac{e(e-1)}{e-g}+\frac{1+1/8}{\m (e-g)}.
\end{alignat*}
We claim that the second derivative $\IG''(g)$ is positive
for $1\leq g\leq e/2$. One finds
\begin{alignat*}1
\IG''(g)=&\frac{2(e^2/g^3+1)(e-g)+2(2g+1-e^2/g^2)}{(e-g)^2}+
\frac{2e(e-1)}{(e-g)^3}\\
\nonumber+&\frac{2(g^2+g+e^2/g+e)}{(e-g)^3}+
\frac{2(1+1/8)}{\m (e-g)^3}.
\end{alignat*}
For $1\leq g\leq e/2$ the last three fractions are certainly
positive and so we may focus on the numerator of the
first fraction. Now if $2g+1-e^2/g^2\geq 0$ the claim
follows at once. If $2g+1-e^2/g^2<0$ it suffices to show
that 
\begin{alignat*}1
(e^2/g^3+1)(e-g)\geq e^2/g^2-2g-1.
\end{alignat*}
With $u=e/g$ the latter is equivalent to
$u^3-u^2+e-g\geq u^2-2g-1$ and this is equivalent
to $u^2(u-2)+e+g+1\geq 0$, which is certainly true
since $1\leq g\leq e/2$ and therefore $2\leq u \leq e$.\\

Thus we have shown that $\IG''(g)>0$ for $1\leq g\leq e/2$
so that $\IG$ is here concave. It suffices to prove $\IG(g)-(n+1)\leq -1/(4\m e)$
for $g=1$, $g=e/2$.
First we use a simple arithmetic argument. Since $n$ is an integer and $n>E=5e/2+4+2/(\m e)$ with denominator dividing
$2\m e$ we see that
\begin{alignat}1
\label{nest}
n+1\geq E+1+1/(2\m e).
\end{alignat}
Now $\IG(e/2)=5e/2+5+2/(\m e)+1/(4\m e)=E+1+1/(4\m e)$ and
thus
\begin{alignat*}1
\IG(e/2)-(n+1)\leq 1/(4\m e)-1/(2\m e)=-1/(4\m e).
\end{alignat*}
Finally 
\begin{alignat*}1
\IG(1)=2e+2+4/(e-1)+9/(8\m (e-1)).
\end{alignat*}
Using (\ref{nest}) again yields
\begin{alignat}1
\label{IG1est}
\IG(1)-(n+1)\leq \frac{4}{e-1}+\frac{9}{8\m(e-1)}-\frac{e}{2}-3-
\frac{2}{\m e}-\frac{1}{2\m e}.
\end{alignat}
First suppose $e=2$. Then (\ref{IG1est}) says
$\IG(1)-(n+1)\leq-1/(8\m)=-1/(4e\m)$. Next suppose $e>2$. Then
the right-hand side in (\ref{IG1est}) is $\leq 4/2+9/(16\m)-e/2-3-5/(2\m e)<
-5/(2\m e)<-1/(4e\m)$.
This completes the proof of the lemma.
\end{rproof}

To show convergence for $\sum_K \Mainterm$ we may use
similar arguments but here we use only $\delta=\delta_1$ instead of $\delta_g$. Let $d=\m e$ so that $[K:\IQ]=d$.
To estimate $V_{\en_K}$ in (\ref{mainterm1}) recall that $V_{\en_K}=V_{\en_K}^{inf}V_{\en_K}^{fin}$.
By (\ref{Vinfbou1}) we have 
\begin{alignat*}1
V_{\en_K}^{inf}\ll (C_{\en_K}^{inf})^{d(n+1)}. 
\end{alignat*}
To estimate $V_{\en_K}^{fin}$ we define the non-zero ideal $\C_0$ by
\begin{alignat}3
\label{C0}
\C_0=\prod_{v\nmid \infty}\pw_v^{-\frac{d_v\log c_v}{\log N\pw_v}}
\end{alignat}
with $c_v$ as in (\ref{defcfin}).
Thus $|\C_0|_v=c_v$ and
\begin{alignat}3
\label{NC0}
N\C_0=(C_{\en_K}^{fin})^d.
\end{alignat}
Let $\D\neq 0$ be a fractional ideal.
Clearly $|\alpha|_v\leq |\C_0^{-1}\D|_v$ for all
non-archimedean $v$ is equivalent to $\alpha \in \C_0^{-1}\D$.
By (\ref{Nineq1}) we conclude
\begin{alignat}3
\label{OG}
\LamenK(\D) \subseteq \sigma(\C_0^{-1}\D)^{n+1}
\end{alignat}
(where $\sigma$ is given by (\ref{sigD})) and therefore
\begin{alignat*}1
\det {\LamenK}(\D)\geq \det \sigma(\C_0^{-1}\D)^{n+1}.
\end{alignat*}
It is well-known (see \cite{23} p.33 (5.2) Satz) that
\begin{alignat*}1
\det \sigma(\C_0^{-1}\D)=
2^{-s_K}\sqrt{|\Delta_K|}N(\D)N(\C_0)^{-1},
\end{alignat*}
where $s_K$ is the number of pairs of complex conjugate embeddings of $K$.
Combining the latter with (\ref{NC0}) we see that
\begin{alignat*}1
\det \sigma(\C_0^{-1}\D)^{n+1}=2^{-s_K(n+1)}|\Delta_K|^{(n+1)/2}N\D^{n+1}(C_{\en_K}^{fin})^{-d(n+1)}.
\end{alignat*}
Inserting
the latter in definition (\ref{defVfin1}) yields
\begin{alignat*}1
V_{\en_K}^{fin}\ll (C_{\en_K}^{fin})^{d(n+1)}.
\end{alignat*}
Now on recalling that $C_{\en_K}=C_{\en_K}^{inf}C_{\en_K}^{fin}$ 
and using (\ref{CenCenk}) we conclude
\begin{alignat*}1
V_{\en_K}\ll C_{\en_K}^{d(n+1)}\leq C_{\en}^{d(n+1)}.
\end{alignat*}
The number of roots of unity $\wK$ in (\ref{Schanuelconst}) is at least $2$.
Furthermore $\zeta_K(n+1)>1$. Hence 
$S_K(n)\ll R_Kh_K|\Delta_K|^{-\frac{n+1}{2}}$.
This together with the above estimate for $V_{\en_K}$
implies $\Mainterm\ll C_{\en}^{d(n+1)}R_Kh_K|\Delta_K|^{-\frac{n+1}{2}}$ 
and since by Siegel-Brauer $R_Kh_K\ll_{\epsilon} |\Delta_K|^{\frac{1}{2}+\epsilon}$ for any positive $\epsilon$
we get 
\begin{alignat}1
\label{Maintermdelta}
\Mainterm \ll_{\epsilon} C_{\en}^{d(n+1)}|\Delta_K|^{-\frac{n}{2}+\epsilon}.
\end{alignat}
\begin{remark}
Let $N_{\Delta}(T)=|\{K\in \coll_e;|\Delta_K|\leq T\}|$. 
Schmidt \cite{55} showed
\begin{alignat}1
\label{Schmidtdisc}
N_{\Delta}(T)\ll T^{\frac{e+2}{4}}.
\end{alignat}
Thus we could apply the dyadic summation lemma
with $\f=|\Delta_K|$ and $b=(e+2)/4$ to see that $\sum_K \Mainterm$
converges for $n>e/2+1$.
\end{remark}
Instead of using Schmidt's bound (\ref{Schmidtdisc}) we will prove a lower bound for $|\Delta_K|$ in terms of $\delta(K/k)$
which might be of interest for its own sake. Then we can apply Lemma \ref{lemmagamma} instead of (\ref{Schmidtdisc}).
\begin{lemma}\label{lemmade4}
We have
\begin{alignat}1
\label{de4}
\delta(K/k)\leq \delta(K/\IQ)\ll |\Delta_K|^{\frac{1}{d}}.
\end{alignat}
\end{lemma}
\begin{rproof}
The lemma is trivially true for $K=k=\IQ$.
However we have by assumption $e\geq 2$
and so $[K:\IQ]=e\m\geq 2$.
The first inequality follows immediately from the definition.
Let $\sigma$ be as in (\ref{sigd}) and suppose
$\alpha$ is a non-zero integer of $K$.
One gets 
\begin{alignat}3\label{heightboundalgint}
\nonumber H(1,\alpha)&=\prod_{i=1}^{d}\max\{1,|\sigma_i(\alpha)|\}^{1/d}\\
\nonumber &\leq \max\{1,\max_{1\leq i\leq d}\{|\sigma_i(\alpha)|\}\}\\
&\leq|\sigma(\alpha)|
\end{alignat}
because $\prod_{1\leq i \leq d}|\sigma_i(\alpha)|\geq 1$.
Let $v_1=\sigma(\alpha_1),...,v_d=\sigma(\alpha_d)$ be linearly independent vectors of the lattice $\sigma \Oseen_K$ with
$|v_i|=\lambda_i$ for the successive minima 
$\lambda_i$ ($i=1,...,[K:\IQ]=d$).
Let us temporarily denote by $b$ the maximum of the degrees of the proper
subfields of $K$. Therefore
$K=\IQ(\alpha_1,...,\alpha_{b+1})$. Next we need to construct
a primitive element in $\Oseen_K$ with small height. 
A standard argument (see \cite{AlgLa} p.244 or Lemma 3.3 in \cite{art1}) yields a primitive
$\alpha=\sum_{j=1}^{b+1}m_j\alpha_j$ with rational integers $0\leq m_j<e$.
Hence by (\ref{heightboundalgint}) we get 
\begin{alignat*}1
H(1,\alpha)\leq |\sigma(\sum_{j=1}^{b+1}m_j\alpha_j)|\leq 
\sum_{j=1}^{b+1}m_j|\sigma(\alpha_j)|\ll \lambda_{b+1}.
\end{alignat*}
We shall estimate $\lambda_{b+1}$:
\begin{alignat*}2
\lambda_{b+1} &= 
\left(\frac{\lambda_1...\lambda_{b}\lambda_{b+1}^{d-b}}
{\lambda_1...\lambda_{b}}\right)^{\frac{1}{d-b}}\\
&\leq  
\left(\frac{\lambda_1...\lambda_d}
{\lambda_1...\lambda_{b}}\right)^{\frac{1}{d-b}}\\
&\ll
\left(\frac{\det(\sigma\Oseen_K)}
{\lambda_1...\lambda_{b}}\right)^{\frac{1}{d-b}}
\\
&=
\left(\frac{ |\Delta_K|^{\frac{1}{2}}}
{2^{s_K}\lambda_1...\lambda_{b}}\right)^{\frac{1}{d-b}}
\\
&\ll |\Delta_K|^{\frac{1}{2(d-b)}}
\end{alignat*}
where we used that $\lambda_1=|\sigma(\alpha_1)|\geq H(1,\alpha_1)\geq 1$.
So all this together implies 
\begin{alignat}1
\label{de3}
\delta(K/\IQ)\ll |\Delta_K|^{\frac{1}{2(d-b)}}.
\end{alignat}
Now $b$ is the degree of a proper subfield. Thus 
$b\leq d/2$ and we get (\ref{de4}).
\end{rproof}
Using Lemma \ref{lemmade4} and (\ref{Maintermdelta}) with
$\epsilon$ replaced by $\epsilon/d$ we deduce  
\begin{alignat}1
\label{deltaDelta1}
\Mainterm &\ll_{\epsilon} C_{\en}^{d(n+1)}\delta(K/k)^{-\frac{dn}{2}+\epsilon}
\end{alignat}
for any positive $\epsilon$.
Choosing $\epsilon=1/2$ we get
\begin{alignat}1
\nonumber \Mainterm &\ll C_{\en}^{d(n+1)}\delta(K/k)^{-\frac{dn}{2}+\frac{1}{2}}.
\end{alignat}
Applying the dyadic summation lemma with
$\f=\delta$ and $b=\gamma$ from Lemma \ref{lemmagamma} we
conclude
\begin{alignat*}2
\sum_{K\in \coll_e} \Mainterm &\ll
C_{\en}^{d(n+1)} \sum_{K\in \coll_e}\delta(K/k)^{-\frac{dn}{2}+\frac{1}{2}}\\
&\ll C_{\en}^{d(n+1)}\sum_{i=1}^{\infty}2^{(-\frac{dn}{2}+\frac{1}{2}+\gamma)i}\\
&\ll C_{\en}^{d(n+1)}
\end{alignat*}
provided $-\frac{dn}{2}+\frac{1}{2}+\gamma<0$, which is 
equivalent to
$n>2e+2+1/d$. But the latter holds since $n>5e/2+4+2/(\m e)$.
This completes the proof of Theorem \ref{main theorem}.

\section{Counting number fields}\label{countingfields}
Using results of the previous section we give simple lower bounds for the growth rate of $N_{\delta}(T)$ and $N_{\Delta}(T)$, the number of field extensions $K/k$ of degree $e$ with $\delta(K/k)\leq T$ or $|\Delta_K|\leq T$. The following corollary shows that the estimates for $N_{\delta}(T)$ are more precise than those available for $N_{\Delta}(T)$. Recall that $e>1$.
\begin{korollar}\label{corlowerdeltabound}
With $\cS=\cS(k,e,1), \CS=\CS(k,e,1)$ and $X_0(k,e,1)$
from (\ref{ThSchm1'}) set
\begin{alignat*}1
c_{\delta}=2^{-5e\m-22}\cS,\text{ } C_{\delta}=\CS \text{ and } T_0=X_0(k,e,1).
\end{alignat*}
Then we have
\begin{alignat*}1
c_{\delta}T^{\m e(e-1)}\leq N_{\delta}(T)\leq C_{\delta}T^{\m e(e+1)}
\end{alignat*}
where the upper bounds holds for $T>0$ and the lower
bound holds for $T\geq T_0$.
\end{korollar}
\begin{rproof}
From the definition it is clear that $Z_H(\IP(K/k),T)>0$ if and only if
$\delta(K/k)\leq T$. Therefore we have
\begin{alignat}1
\label{Ndeltaexpr}
N_{\delta}(T)=\sum_{K \in \coll_e \atop \delta(K/k)\leq T}1=\sum_{K \in \coll_e \atop \delta(K/k)\leq T}\frac{Z_{H}(\IP(K/k),T)}{Z_{H}(\IP(K/k),T)}.
\end{alignat}
Using the equivalence above once again, we see that
the term on the far right-hand side of (\ref{Ndeltaexpr}) is
\begin{alignat*}1
&\geq (\sup_{K \in \coll_e}\{Z_{H}(\IP(K/k),T)\})^{-1}\sum_{K \in \coll_e \atop \delta(K/k)\leq T}Z_{H}(\IP(K/k),T)\\
&= (\sup_{K \in \coll_e}\{Z_{H}(\IP(K/k),T)\})^{-1}\sum_{K \in \coll_e}Z_{H}(\IP(K/k),T)\\
&=(\sup_{K \in \coll_e}\{Z_{H}(\IP(K/k),T)\})^{-1}
Z_{H}(\IP(k;e),T).
\end{alignat*}
Now $Z_{H}(\IP(K/k),T)\leq Z_{H}(\IP(K;1),T)$
and by the upper bound in (\ref{ThSchm1'})
and recalling that $[K:\IQ]=e\m$ we get 
\begin{alignat*}1
Z_{H}(\IP(K;1),T)\leq \CS(K,1,1) T^{2\m e}=2^{5e\m+22}T^{2\m e}.
\end{alignat*}
The lower bound in (\ref{ThSchm1'})
with $\cS=\cS(k,e,1)$
yields
\begin{alignat*}1
Z_{H}(\IP(k;e),T)\geq \cS T^{\m e(e+1)}
\end{alignat*}
for $T\geq X_0(k,e,1)=T_0$.
Hence
\begin{alignat*}1
N_{\delta}(T)\geq (2^{5e\m+22}T^{2\m e})^{-1}\cS T^{\m e(e+1)}=c_{\delta}T^{\m e(e-1)}
\end{alignat*}
for $T\geq T_0$. On the other hand Lemma \ref{lemmagamma} tells us that
\begin{alignat*}1
N_{\delta}(T)\leq C_{\delta} T^{\m e(e+1)}
\end{alignat*}
for $T>0$.
\end{rproof}
Corollary \ref{corlowerdeltabound} combined with the lower bound
(\ref{de11}) for $\delta$ in terms of $|\Delta_K|$ yields
\begin{korollar}\label{corlowerDeltabound}
There are positive constants $\cDe=\cDe(k,e)$ and $T_1=T_1(k,e)$
depending only on $k,e$ such that
\begin{alignat*}1
N_{\Delta}(T)\geq \cDe T^{1/2}
\end{alignat*}
for $T\geq T_1$.
\end{korollar}
\begin{rproof}
From (\ref{de11}) we know that there is a positive constant
$\ceins=\ceins(k,e)$ depending only on $k,e>1$ such that
$\delta(K/k)\geq \ceins|\Delta_K|^{1/(2e(e-1)\m)}$.
Using Corollary \ref{corlowerdeltabound} and setting
$\cDe=c_{\delta}\ceins^{\m e(e-1)}$, $T_1=(T_0/\ceins)^{2e(e-1)\m}$
we conclude
\begin{alignat*}1
N_{\Delta}(T)\geq N_{\delta}(\ceins T^{1/(2e(e-1)\m)})
\geq  \cDe T^{1/2}
\end{alignat*}
provided $T\geq T_1$.
\end{rproof}
Ellenberg and Venkatesh's Theorem 1.1 in \cite{56}
shows that the exponent $1/2$ in Corollary \ref{corlowerDeltabound}
can be replaced by $1/2+1/e^2$ and according to Linnik's Conjecture (see \cite{56} p.723) the correct exponent is $1$. Although the general Linnik Conjecture is known 
to be true only for $e\leq 3$ 
the exponent $1/2$ can always be increased to $1$  if  $e$ is even
or a multiple of $3$ (see \cite{56} pp. 723,724).\\

What about upper bounds for $N_{\Delta}(T)$?
From (\ref{de4}) we know that there is a positive
constant $\czwei=\czwei(d)$ depending only on $d=e\m$
such that
\begin{alignat*}1
\delta(K/k)\leq
\czwei|\Delta_K|^{\frac{1}{e\m}}.
\end{alignat*}
Thus we get
\begin{alignat*}1
N_{\Delta}(T)\leq N_{\delta}(\czwei T^{\frac{1}{e\m}})
\leq  C_{\delta}\czwei^{\m e(e+1)}T^{e+1}
\end{alignat*}
for $T>0$.
But Schmidt's bound (\ref{Schmidtdisc}) has the much better exponent $(e+2)/4$ on $T$.

\bibliographystyle{siam}
\bibliography{literature}

\end{document}